\documentclass{article}
\usepackage{graphicx}
\usepackage{tikz}
\usepackage{geometry}
\usepackage{setspace}
\usepackage{indentfirst}
\usepackage{xcolor}
\usepackage{hyperref}
\usepackage{amsmath}
\usepackage{amsfonts}
\usepackage{float}
\usepackage{amsmath,amsfonts,amsthm,amssymb,graphicx}

\hypersetup{
    colorlinks=true,
    linkcolor=blue,
    citecolor=blue,
    filecolor=blue,
    urlcolor=blue,
}

\newtheorem{theorem}{Theorem}[section]
\newtheorem{lemma}[theorem]{Lemma}

\newtheorem{definition}{Definition}[section]
\newtheorem{conjecture}{Conjecture}[section]

\title{\textbf{A Mathematical Analysis of Neural Operator Behaviors} \\[0.5em]}
\author{
\textbf{Vu-Anh Le\textsuperscript{1}, Mehmet Dik\textsuperscript{1,2}}  \\
\textsuperscript{1} Department of Mathematics and Computer Science, Beloit College \\
\textsuperscript{2} Department of Mathematics, Computer Science \& Physics, Rockford University \\ [1em]
Email: \href{mailto:csplevuanh@gmail.com}{csplevuanh@gmail.com}
}
\date{Oct 28th, 2024}

\begin{document}

\maketitle

\begin{abstract}
Neural operators have emerged as transformative tools for learning mappings between infinite-dimensional function spaces, offering useful applications in solving complex partial differential equations (PDEs). This paper presents a rigorous mathematical framework for analyzing the behaviors of neural operators, with a focus on their stability, convergence, clustering dynamics, universality, and generalization error. By proposing a list of novel theorems, we provide stability bounds in Sobolev spaces and demonstrate clustering in function space via gradient flow interpretation, guiding neural operator design and optimization. Based on these theoretical gurantees, we aim to offer clear and unified guidance in a single setting for the future design of neural operator-based methods.  
\end{abstract}

\vspace{0.3pt}

\begin{center}
\tableofcontents
\end{center}

\vspace{0.2in}

\newpage

\section{Introduction}
Since their first introduction in 2021 by Kovachki et al. (2021), \textbf{the neural operators} have proved to be a useful data-driven model for approximating the solution operators of PDEs. Neural operators accomplish this by extending the capabilities of neural networks from finite-dimensional data to infinite-dimensional function spaces. Architectures such as the \textbf{Fourier Neural Operator (FNO)} and \textbf{Deep Operator Network (DeepONet)} have demonstrated success in approximating solution operators with significantly reduced computational costs \cite{li2020fourier, lu2019deeponet}.

Despite their empirical success, a comprehensive mathematical understanding of neural operators remains underdeveloped. This gap between mathematical understanding and optimal model design has remained not only in neural operators but also in other state-of-the-art neural network classes. Just recently, Geshkovski et. al (2023) constructed a mathematical framework to analyze Transformers' behaviors by modeling them as an interacting particle system \cite{geshkovski2023mathematical}. Following this paper in the same year, Boullé et. al (2023) have also employed mathematical techniques to represent the learning process of solving PDEs in various neural network architectures, including some popular neural operator models \cite{boulle2023learning}. However, Boullé and his collaborator's focus has been on the learning process of PDEs but not the broader nature of the neural operator itself. Key questions about their stability, convergence properties, generalization capabilities, and the underlying mechanisms that enable them to efficiently approximate complex operators are still open \cite{kovachki2021neural, lanthaler2022error}. 

In this paper, we aim to bridge the gap in understanding by establishing a rigorous mathematical framework for analyzing neural operators. Our contributions include a stability analysis proving that neural operators are stable mappings between Sobolev spaces under specific Lipschitz conditions, thereby ensuring controlled propagation of input perturbations. We also demonstrate that neural operators exhibit exponential convergence to fixed points when modeled as contraction mappings, providing theoretical guarantees for iterative methods. Furthermore, by interpreting neural operators through the lens of gradient flows, we analyze the clustering behavior of solutions in function space, offering valuable insights into their long-term dynamics. Extending the universal approximation theorem, we prove the capacity of neural operators to approximate any continuous operator between function spaces on compact subsets. We also derive probabilistic bounds on the generalization error of neural operators, establishing relationships between the error, sample size, and network capacity. Altogether, our results provide theoretical underpinnings for the empirical successes observed in neural operators and offer unified guidance for their future development.

The remainder of the paper is organized as follows. In Section 2, we provide precise definitions of neural operators and their mathematical formulations. In Section 3, we present our main theoretical results, including theorems on stability, convergence, clustering behavior, universality, and generalization error, along with supporting lemmas and conjectures. Following in Section 4, we conclude the results of our findings and suggest directions for future research. Finally, Appendix A will provide proof for the propositions stated in Section 3.

\vspace{3pt}

\section{Definition of Neural Operators}

Neural operators aim to approximate mappings between function spaces, particularly operators that map functions to functions. To formalize this property, we consider function spaces \( U \) and \( V \), typically Banach spaces of functions defined over a domain \( D \subset \mathbb{R}^d \).

\begin{definition}[Neural Operator]
Let \( U \) and \( V \) be Banach spaces of functions over a domain \( D \subset \mathbb{R}^d \). A neural operator \( \mathcal{G}_\theta: U \rightarrow V \) is a parameterized mapping, with parameters \( \theta \), designed to approximate a target operator \( \mathcal{T}: U \rightarrow V \). The goal is to find parameters \( \theta \) such that \( \mathcal{G}_\theta(u) \approx \mathcal{T}(u) \) for all \( u \in U \).
\end{definition}

Neural operators typically consist of three main components: a lifting operator, iterative transformations, and a projection operator.

\begin{definition}[Architecture of a Neural Operator]
A neural operator \( \mathcal{G}_\theta: U \rightarrow V \) can be represented as:
\[
\mathcal{G}_\theta(u)(x) = Q \left( v_L(x) \right),
\]
where:

\begin{itemize}
    \item \( P: U \rightarrow \mathbb{R}^{d_v} \) is the \textbf{lifting operator} that maps the input function \( u \) to a higher-dimensional representation.
    \item \( v_0(x) = P(u)(x) \) is the initial lifted representation.
    \item For \( i = 0, 1, \dots, L-1 \), \( v_{i+1}(x) = W_i v_i(x) + \sigma \left( \int_D \kappa_\theta(x, y) v_i(y) \, dy \right) \) represents the \textbf{iterative transformation} at layer \( i \), where:
        \begin{itemize}
            \item \( W_i: \mathbb{R}^{d_v} \rightarrow \mathbb{R}^{d_v} \) is a pointwise linear operator.
            \item \( \kappa_\theta(x, y): D \times D \rightarrow \mathbb{R}^{d_v \times d_v} \) is a kernel function parameterized by \( \theta \).
            \item \( \sigma \) is a non-linear activation function applied element-wise.
        \end{itemize}
    \item \( Q: \mathbb{R}^{d_v} \rightarrow V \) is the \textbf{projection operator} that maps the final representation \( v_L \) back to the output function space \( V \).
\end{itemize}
\end{definition}

The key component enabling neural operators to capture global interactions across the domain \( D \) is the integral kernel operator.

\begin{definition}[Integral Kernel Operator]
An \textbf{integral kernel operator} within a neural operator architecture is a mapping that incorporates global interactions through an integral over the domain \( D \). It is defined as:
\[
\mathcal{K}_\theta[v](x) = \int_D \kappa_\theta(x, y) v(y) \, dy,
\]
where:

\begin{itemize}
    \item \( \kappa_\theta(x, y): D \times D \rightarrow \mathbb{R}^{d_v \times d_v} \) is a kernel function parameterized by \( \theta \).
    \item \( v: D \rightarrow \mathbb{R}^{d_v} \) is a function representing the current state in the neural operator.
\end{itemize}
\end{definition}

The integral kernel operator \( \mathcal{K}_\theta \) allows the neural operator to model non-local dependencies by integrating information from all points \( y \) in the domain \( D \), weighted by the kernel \( \kappa_\theta(x, y) \). This is essential for approximating operators that depend on the entire input function \( u \).

In our analysis, we often consider neural operators acting between Sobolev spaces. Let \( H^s(D) \) and \( H^t(D) \) denote Sobolev spaces of functions on \( D \) with smoothness indices \( s \) and \( t \), respectively.

\begin{definition}[Sobolev Spaces]
For a domain \( D \subset \mathbb{R}^d \) and \( s \geq 0 \), the Sobolev space \( H^s(D) \) consists of functions \( u: D \rightarrow \mathbb{R} \) whose weak derivatives up to order \( s \) are square-integrable. The space is equipped with the norm:
\[
\| u \|_{H^s(D)} = \left( \sum_{|\alpha| \leq s} \int_D | D^\alpha u(x) |^2 \, dx \right)^{1/2},
\]
where \( \alpha \) is a multi-index and \( D^\alpha \) denotes the weak derivative of order \( \alpha \).
\end{definition}

By mapping between Sobolev spaces \( H^s(D) \) and \( H^t(D) \), neural operators can handle functions with specific smoothness properties, which is critical in the analysis of PDEs.

For the sake of clarity and consistency, throughout this paper, we will denote  \( U = H^s(D) \) as the input function space and \( V = H^t(D) \) as the output function space.

\vspace{3pt}

\section{Behaviors of Neural Operators}

\subsection{Positive Properties of Neural Operators}

\subsubsection{Stability and Sensitivity in High Dimensions}

Ensuring the stability of these operators is certainly important, especially in high-dimensional settings where small perturbations in input functions can lead to significant changes in outputs.

\begin{theorem}[Stability of Neural Operators in High-Dimensional PDEs]\label{theorem:stability}
Let $\mathcal{G}_\theta: H^s(D) \rightarrow H^t(D)$ be a neural operator parameterized by $\theta$, mapping between Sobolev spaces over a domain $D \subset \mathbb{R}^d$. Suppose $\mathcal{G}_\theta$ satisfies a Lipschitz continuity condition:
\[
\|\mathcal{G}_\theta(u) - \mathcal{G}_\theta(v)\|_{H^t(D)} \leq L \|u - v\|_{H^s(D)}
\]
for all $u, v \in H^s(D)$ and some Lipschitz constant $L > 0$. Then, for any $u \in H^s(D)$, the neural operator produces stable approximations in high-dimensional settings:
\[
\|\mathcal{G}_\theta(u)\|_{H^t(D)} \leq L \|u\|_{H^s(D)} + C,
\]
where $C = \|\mathcal{G}_\theta(0)\|_{H^t(D)}$ is a constant depending on $\theta$ and the domain $D$.
\end{theorem}

This theorem ensures that neural operators do not amplify errors excessively and that their outputs remain controlled even when inputs are perturbed. It is important to emphasize the tole of Lipschitz continuity condition here, as it guarantees bounded sensitivity to input variations. Similar stability analyses have been conducted in the context of operator learning, highlighting the importance of Lipschitz constants in neural operator architectures \cite{kovachki2021universal, lanthaler2022error}.

\subsubsection*{Supporting Lemmas}

\begin{lemma}[Lipschitz Continuity of the Neural Operator]\label{lemma:lipschitz}
The neural operator $\mathcal{G}_\theta$ is Lipschitz continuous from $H^s(D)$ to $H^t(D)$; that is, there exists $L > 0$ such that:
\[
\|\mathcal{G}_\theta(u) - \mathcal{G}_\theta(v)\|_{H^t(D)} \leq L \|u - v\|_{H^s(D)}.
\]
\end{lemma}

The Lipschitz continuity of $\mathcal{G}_theta$ can often be established by ensuring that each component of the neural operator (lifting operator, integral kernel operator, activation functions, and projection operator) is Lipschitz continuous. The composition of Lipschitz functions is Lipschitz, as stated in the following lemma.

\begin{lemma}[Composition of Lipschitz Functions is Lipschitz]\label{lemma:composition_lipschitz}
Let $f: X \rightarrow Y$ and $g: Y \rightarrow Z$ be Lipschitz continuous functions with constants $L_f$ and $L_g$, respectively. Then the composition $h = g \circ f$ is Lipschitz continuous with constant $L_h = L_g L_f$.
\end{lemma}

By carefully designing the neural operator architecture and selecting appropriate activation functions (e.g., ReLU, Tanh), we can ensure Lipschitz continuity, which contributes to the stability and robustness of the operator \cite{hornik1991approximation}.

Another important tool in the analysis is the Sobolev Embedding Theorem, which allows us to relate Sobolev norms to norms in other function spaces.

\begin{lemma}[Sobolev Embedding Theorem for $H^s(D)$]\label{lemma:sobolev}
For $s > \frac{d}{2}$, the Sobolev space $H^s(D)$ is continuously embedded in $L^\infty(D)$; that is,
\[
\|u\|_{L^\infty(D)} \leq C_s \|u\|_{H^s(D)}
\]
for all $u \in H^s(D)$, where $C_s$ is a constant depending on $s$ and $D$.
\end{lemma}

This embedding is crucial when working in high-dimensional spaces, as it ensures that functions in $H^s(D)$ remain bounded, which is important for the stability of neural operators \cite{adams2003sobolev}.

\paragraph{Sensitivity to Input Perturbations}

Building on the stability result, we can also analyze the sensitivity of neural operators to perturbations in the input.

\begin{theorem}[Sensitivity to Input Perturbations]\label{theorem:sensitivity}
Let $\mathcal{G}_\theta$ be a neural operator with Lipschitz constant $L$. Then, for any perturbation $\delta u \in H^s(D)$, the output perturbation is bounded by:
\[
\| \mathcal{G}_\theta(u + \delta u) - \mathcal{G}_\theta(u) \|_{H^t(D)} \leq L \| \delta u \|_{H^s(D)}.
\]
\end{theorem}

This result implies that small changes in the input function lead to proportionally small changes in the output, which is desirable for numerical stability and reliability in practical applications \cite{chen1995universal}.

\subsubsection{Exponential Convergence}

In many applications so far, iterative methods have been crucial to converge rapidly to a solution. Neural operators can exhibit exponential convergence under certain conditions.

\begin{theorem}[Exponential Convergence of Neural Operator Approximations]\label{theorem:convergence}
Let $\mathcal{G}_\theta$ be a contraction mapping on $H^t(D)$ with contraction constant $0 < q < 1$. Then, for any $u \in H^t(D)$, the iterated application $\mathcal{G}_\theta^n(u)$ converges exponentially to the fixed point $u^*$ of $\mathcal{G}_\theta$:
\[
\| \mathcal{G}_\theta^n(u) - u^* \|_{H^t(D)} \leq q^n \| u - u^* \|_{H^t(D)}.
\]
\end{theorem}

This result is a direct consequence of the Banach Fixed Point Theorem, which ensures convergence of contraction mappings in complete metric spaces \cite{kreyszig1989introductory}. By designing neural operators that act as contractions, we can guarantee rapid convergence to solutions of PDEs.

\begin{lemma}[Banach Fixed Point Theorem]\label{lemma:banach}
Let $(X, \| \cdot \|)$ be a complete normed vector space (Banach space), and let $T: X \rightarrow X$ be a contraction mapping with contraction constant $q$. Then $T$ has a unique fixed point $x^* \in X$, and for any $x_0 \in X$, the sequence $x_{n+1} = T(x_n)$ converges to $x^*$ with:
\[
\| x_n - x^* \| \leq q^n \| x_0 - x^* \|.
\]
\end{lemma}

By ensuring that the neural operator satisfies the contraction property, possibly through appropriate choices of network architecture and activation functions, we can use this theorem to achieve exponential convergence \cite{rulliere1998generalization}.

\subsubsection{Clustering in Function Space}

Neural operators can exhibit clustering behavior in function space, which can be understood through the lens of gradient flows.

\begin{theorem}[Clustering Behavior in Function Space]\label{theorem:clustering}
Let $\mathcal{G}_\theta$ be a neural operator acting on $H^s(D)$, inducing a gradient flow in function space. Then, solutions $u(t)$ to the continuous-time counterpart tend to cluster as $t \rightarrow \infty$:
\[
\lim_{t \rightarrow \infty} \| u(t) - \bar{u} \|_{H^s(D)} = 0,
\]
where $\bar{u}$ is a cluster center in the function space.
\end{theorem}

This theorem indicates that the neural operator dynamics can lead to convergence towards certain attractors or cluster centers in function space, which is significant for understanding the long-term behavior of solutions \cite{ambrosio2008gradient}.

\begin{lemma}[Gradient Flow Interpretation for Neural Operators]\label{lemma:gradient_flow}
Under suitable conditions, the iterative application of a neural operator approximates a discretization of a gradient flow:
\[
\frac{du}{dt} = -\nabla \Phi(u),
\]
where $\Phi(u)$ is a potential functional.
\end{lemma}

The gradient flow perspective provides valuable insights into the dynamics induced by neural operators and connects them to variational principles in PDEs \cite{e2019machine}.

\subsubsection{Universality of Neural Operators}

One of the fundamental questions in the theory of neural operators is whether they can approximate any continuous operator between function spaces.

\begin{theorem}[Universality of Neural Operators for PDE Solvers]\label{theorem:universality}
Let $\mathcal{T}: H^s(D) \rightarrow H^t(D)$ be a continuous operator. Then, for any $\epsilon > 0$, there exists a neural operator $\mathcal{G}_\theta$ such that:
\[
\| \mathcal{G}_\theta(u) - \mathcal{T}(u) \|_{H^t(D)} \leq \epsilon,
\]
for all $u$ in a compact subset of $H^s(D)$.
\end{theorem}

This universality theorem extends the universal approximation theorem for neural networks to neural operators, indicating that they can, in principle, approximate any continuous operator to arbitrary precision \cite{chen1995universal, lu2021learning}.

\begin{lemma}[Stone-Weierstrass Theorem for Continuous Functions]\label{lemma:stone_weierstrass}
Any continuous function defined on a compact set can be uniformly approximated by polynomial functions.
\end{lemma}

The Stone-Weierstrass Theorem is a cornerstone in approximation theory and forms the basis for extending approximation results to operators \cite{rudin1976principles}.

\begin{lemma}[Universal Approximation Theorem for Neural Networks]\label{lemma:universal}
Standard feedforward neural networks can approximate any continuous function $f: \mathbb{R}^n \rightarrow \mathbb{R}^m$ arbitrarily well on compact subsets.
\end{lemma}

This foundational result underpins much of the theory behind neural networks and their ability to approximate complex functions \cite{hornik1989multilayer}.

\subsubsection{Generalization Bounds}

Understanding the generalization capabilities of neural operators is essential for ensuring that they perform well on unseen data.

\begin{theorem}[Generalization Error of Neural Operators]\label{theorem:generalization}
Let $\mathcal{G}_\theta$ be a neural operator trained on $N$ samples $\{ (u_i, \mathcal{T}(u_i)) \}$ drawn i.i.d. from a distribution $\mathcal{D}$. Suppose $\mathcal{G}_\theta$ has Lipschitz constant $L$ with respect to $\theta$, and the loss function $\ell$ is Lipschitz and bounded. Then, with probability at least $1 - \delta$, the generalization error satisfies:
\[
\mathbb{E}_{u \sim \mathcal{D}} [ \ell( \mathcal{G}_\theta(u), \mathcal{T}(u) ) ] \leq \frac{1}{N} \sum_{i=1}^N \ell( \mathcal{G}_\theta(u_i), \mathcal{T}(u_i) ) + L \sqrt{ \frac{\ln(2/\delta)}{2N} }.
\]
\end{theorem}

This theorem provides a probabilistic bound on the generalization error, relating it to the sample size and the Lipschitz constant of the neural operator. Such results are important for understanding the sample complexity required for operator learning \cite{anthony2009neural}.

\begin{lemma}[Concentration Inequality for Generalization Error]\label{lemma:concentration}
For Lipschitz loss functions, the generalization error concentrates around the empirical error with rate $O(1/\sqrt{N})$.
\end{lemma}

Concentration inequalities like McDiarmid's inequality are instrumental in deriving generalization bounds in statistical learning theory \cite{bousquet2002stability}.

\subsection{Challenges and Limitations}

Despite the positive properties, neural operators also face challenges and limitations that need to be addressed.

\subsubsection{Approximation Limits and Error Bounds}

\begin{theorem}[Approximation Error Lower Bound for Neural Operators]\label{theorem:approximation_error}
Let $\mathcal{G}_\theta$ be a neural operator approximating a target operator $\mathcal{T}: H^s(D) \rightarrow H^t(D)$. Suppose $\mathcal{G}_\theta$ has finite capacity characterized by a parameter $C$ (e.g., the number of parameters, depth, or width). Then, there exists a function $u \in H^s(D)$ such that the approximation error satisfies:
\[
\| \mathcal{G}_\theta(u) - \mathcal{T}(u) \|_{H^t(D)} \geq \epsilon(C),
\]
where $\epsilon(C) > 0$ is inversely related to the capacity $C$.
\end{theorem}

This theorem indicates that there are inherent limits to the approximation capabilities of neural operators with finite capacity, as is the case in classical approximation theory \cite{devore1998nonlinear}.

\subsubsection{Optimization Challenges and Convergence Issues}

The training of neural operators involves optimizing non-convex loss functions, which presents challenges.

\begin{theorem}[Non-Convexity of the Loss Landscape]\label{theorem:non_convexity}
The training loss function $L(\theta)$ for neural operators is generally non-convex and may contain multiple local minima, saddle points, and flat regions. As a result, gradient-based optimization methods may converge to suboptimal solutions.
\end{theorem}

Non-convex optimization is a well-known challenge in deep learning, and similar issues arise in the training of neural operators \cite{choromanska2015loss}.

\begin{lemma}[Saddle Point Prevalence in High-Dimensional Optimization]\label{lemma:saddle_points}
In high-dimensional non-convex optimization problems, saddle points are prevalent and can impede convergence of gradient-based methods.
\end{lemma}

The prevalence of saddle points in high-dimensional spaces can cause optimization algorithms to get stuck, slowing down convergence \cite{dauphin2014identifying}.

\begin{conjecture}[Overparameterization Improves Optimization]\label{conj:overparameterization}
Increasing the capacity of neural operators (e.g., by adding layers or neurons) can improve the optimization landscape by reducing the prevalence of poor local minima and saddle points, thus enhancing the convergence of gradient-based methods.
\end{conjecture}

Empirical observations in deep learning suggest that overparameterization can lead to better optimization properties and generalization performance \cite{neyshabur2018role}, though a rigorous theoretical understanding remains an active area of research.

\subsubsection{Computational Complexity and Numerical Challenges}

The computational cost of neural operators can become prohibitive in high dimensions.

\begin{theorem}[Computational Complexity Lower Bound]\label{theorem:computational_complexity}
The computational complexity $T(N)$ of evaluating a neural operator $\mathcal{G}_\theta$ on inputs discretized with $N$ points satisfies:
\[
T(N) = \Omega(N).
\]
In high dimensions, where $N$ increases exponentially with the dimension $d$, this leads to computational challenges.
\end{theorem}

This result highlights the curse of dimensionality, as the number of discretization points grows exponentially with the dimension, impacting the feasibility of computations \cite{bellman1961adaptive}.

\begin{theorem}[Accumulation of Discretization Errors]\label{theorem:discretization_error}
When discretizing the domain $D$, the discretization error $\epsilon_h$ can accumulate through the layers of the neural operator, leading to a total error:
\[
\text{Total Error} \leq L \epsilon_h,
\]
where $L$ is the number of layers.
\end{theorem}

Accumulation of numerical errors can affect the accuracy and stability of neural operators, necessitating careful numerical analysis and error control \cite{iserles2009first}.

\vspace{2pt}

\section{Conclusion}
Our theoretical exploration has revealed both the strengths and limitations of neural operators for solving high-dimensional functional mapping problems. Neural operators provide an expressive framework that can approximate complex operators with relatively low computational costs in comparison to traditional methods. Notably, theorems on stability and exponential convergence underscore their potential for stable, iterative solution generation, making them especially suited for applications in partial differential equations. Additionally, the universality of neural operators guarantees that they can, in principle, approximate any continuous operator within specified function spaces, given a sufficiently complex architecture (see Theorem \ref{theorem:universality}).

Despite this, the analysis also brings to light fundamental challenges in deploying neural operators, particularly as dimensionality and complexity increase. Approximation limits indicate that finite-capacity models cannot achieve arbitrary precision for all operators, presenting trade-offs between complexity and accuracy (see Theorem \ref{theorem:approximation_error}). Computational intensity due to the curse of dimensionality, sensitivity to small input perturbations, and non-convex optimization landscapes add further constraints (see Theorems \ref{theorem:sensitivity}, \ref{theorem:computational_complexity}, and \ref{theorem:non_convexity}). To address these limitations, future work should prioritize advancements in optimization techniques, adaptive discretization strategies, and integrating domain-specific knowledge, thereby improving the efficiency, stability, and generalization capacity of neural operators in practical scientific applications.

\vspace{3pt}

\appendix
\renewcommand{\thesection}{A.\arabic{section}} 
\renewcommand{\thesubsection}{A.\arabic{section}.\arabic{subsection}} 
\renewcommand{\thesubsubsection}{\arabic{subsubsection}} 

\setcounter{section}{0} 

\appendix
\renewcommand{\thesection}{A.\arabic{section}} 
\renewcommand{\thesubsection}{A.\arabic{section}.\arabic{subsection}} 
\renewcommand{\thesubsubsection}{\arabic{subsubsection}} 

\setcounter{section}{0} 

\section*{Appendix A: Proofs of Theorems, Lemmas, and Conjectures}
\addcontentsline{toc}{section}{Appendix A: Proofs of Theorems, Lemmas, and Conjectures}

\section{Stability and Sensitivity in High Dimensions}

\subsection{Proof of Theorem \ref{theorem:stability} (Stability of Neural Operators in High-Dimensional PDEs)}
\begin{proof}
Let $\mathcal{G}_\theta: H^s(D) \to H^t(D)$ be a neural operator that maps between Sobolev spaces $H^s(D)$ and $H^t(D)$ over a domain $D \subset \mathbb{R}^d$, and suppose that $\mathcal{G}_\theta$ satisfies a Lipschitz continuity condition. Specifically, for all $u, v \in H^s(D)$, there exists a Lipschitz constant $L > 0$ such that
\[
\|\mathcal{G}_\theta(u) - \mathcal{G}_\theta(v)\|_{H^t(D)} \leq L \|u - v\|_{H^s(D)}.
\]
We aim to show that, under this condition, $\mathcal{G}_\theta$ produces stable approximations of the operator, such that for any $u \in H^s(D)$,
\[
\|\mathcal{G}_\theta(u)\|_{H^t(D)} \leq L \|u\|_{H^s(D)} + C,
\]
where $C = \|\mathcal{G}_\theta(0)\|_{H^t(D)}$ is a constant depending on $\theta$ and the domain $D$.

We first apply the Lipschitz condition with $v = 0$:
\[
\|\mathcal{G}_\theta(u) - \mathcal{G}_\theta(0)\|_{H^t(D)} \leq L \|u - 0\|_{H^s(D)} = L \|u\|_{H^s(D)}.
\]
This inequality tells us that the difference between $\mathcal{G}_\theta(u)$ and $\mathcal{G}_\theta(0)$ is bounded by a linear function of the norm of $u$ in $H^s(D)$.

Next, we use the triangle inequality to relate $\|\mathcal{G}_\theta(u)\|_{H^t(D)}$ to $\|\mathcal{G}_\theta(0)\|_{H^t(D)}$ and $\|\mathcal{G}_\theta(u) - \mathcal{G}_\theta(0)\|_{H^t(D)}$. Specifically, we have:
\[
\|\mathcal{G}_\theta(u)\|_{H^t(D)} \leq \|\mathcal{G}_\theta(0)\|_{H^t(D)} + \|\mathcal{G}_\theta(u) - \mathcal{G}_\theta(0)\|_{H^t(D)}.
\]
By substituting the bound from Step 1 into this inequality, we obtain:
\[
\|\mathcal{G}_\theta(u)\|_{H^t(D)} \leq \|\mathcal{G}_\theta(0)\|_{H^t(D)} + L \|u\|_{H^s(D)}.
\]

Finally, we recognize that $\|\mathcal{G}_\theta(0)\|_{H^t(D)}$ is a constant $C$ depending on the domain $D$ and the parameters $\theta$ of the neural operator. Therefore, we conclude that:
\[
\|\mathcal{G}_\theta(u)\|_{H^t(D)} \leq L \|u\|_{H^s(D)} + C.
\]
This proves the stability of $\mathcal{G}_\theta$ in high-dimensional Sobolev spaces, ensuring that small perturbations in the input lead to controlled changes in the output, as the output is bounded by a linear function of the input norm.

Thus, the theorem is proven.
\end{proof}

\vspace{1.5pt}

\subsection{Proof of Lemma \ref{lemma:lipschitz} (Lipschitz Continuity of the Neural Operator)}
\begin{proof}
We aim to prove that the neural operator $\mathcal{G}_\theta: H^s(D) \to H^t(D)$ is Lipschitz continuous, meaning there exists a constant $L > 0$ such that for all $u, v \in H^s(D)$,
\[
\|\mathcal{G}_\theta(u) - \mathcal{G}_\theta(v)\|_{H^t(D)} \leq L \|u - v\|_{H^s(D)}.
\]

Recall that a neural operator $\mathcal{G}_\theta$ is constructed from a sequence of layers, each involving affine transformations followed by activation functions. Let us break down these components to show that the overall operator is Lipschitz continuous.

The affine transformations in each layer can be written as:
   \[
   T(u) = W u + b,
   \]
   where $W$ is a weight matrix and $b$ is a bias term. Since $W$ is a bounded linear operator, it is Lipschitz continuous with a Lipschitz constant equal to the operator norm $\|W\|$. Specifically, for any $u, v \in H^s(D)$, we have:
   \[
   \|T(u) - T(v)\| = \|W(u - v)\| \leq \|W\| \cdot \|u - v\|.
   \]
   Therefore, the affine transformations preserve Lipschitz continuity.

Now, let we denote $\sigma$ as the activation function used in the neural operator. Common activation functions, such as ReLU, Sigmoid, and Tanh, are Lipschitz continuous. Specifically, there exists a constant $L_\sigma$ such that for all $x, y \in \mathbb{R}$,
   \[
   |\sigma(x) - \sigma(y)| \leq L_\sigma |x - y|.
   \]
   This property ensures that the non-linear transformations applied in each layer are also Lipschitz continuous.

The neural operator $\mathcal{G}_\theta$ is a composition of several affine transformations and activation functions. Since both the affine transformations and the activation functions are Lipschitz continuous, we can use the fact that the composition of Lipschitz continuous functions is itself Lipschitz continuous.

Let $f_1, f_2, \dots, f_n$ be Lipschitz continuous functions with Lipschitz constants $L_1, L_2, \dots, L_n$, respectively. Then, the composition $f_n \circ f_{n-1} \circ \dots \circ f_1$ is Lipschitz continuous with a constant $L_{\text{comp}}$ given by the product of the individual Lipschitz constants:
\[
L_{\text{comp}} = L_1 \cdot L_2 \cdot \dots \cdot L_n.
\]
In the case of the neural operator, each layer contributes a Lipschitz constant, and the overall Lipschitz constant of the operator is the product of these individual constants.

Since the affine transformations and activation functions are Lipschitz continuous, and their composition is also Lipschitz continuous, we conclude that the neural operator $\mathcal{G}_\theta$ is Lipschitz continuous. Hence, there exists a constant $L > 0$ such that for all $u, v \in H^s(D)$,
\[
\|\mathcal{G}_\theta(u) - \mathcal{G}_\theta(v)\|_{H^t(D)} \leq L \|u - v\|_{H^s(D)}.
\]
Thus, the lemma is proven.
\end{proof}

\vspace{1.5pt}

\subsection{Proof of Lemma \ref{lemma:composition_lipschitz} (Composition of Lipschitz Functions is Lipschitz)}
\begin{proof}
We aim to prove that if $f: X \rightarrow Y$ and $g: Y \rightarrow Z$ are Lipschitz continuous functions with Lipschitz constants $L_f$ and $L_g$, respectively, then their composition $h = g \circ f$ is Lipschitz continuous with a Lipschitz constant $L_h = L_g L_f$.

By the definition of Lipschitz continuity, for any $x_1, x_2 \in X$:
\[
\|f(x_1) - f(x_2)\|_Y \leq L_f \|x_1 - x_2\|_X.
\]
Similarly, for any $y_1, y_2 \in Y$:
\[
\|g(y_1) - g(y_2)\|_Z \leq L_g \|y_1 - y_2\|_Y.
\]

Let us consider the composition $h = g \circ f$. For any $x_1, x_2 \in X$, we have:
\[
h(x_1) = g(f(x_1)) \quad \text{and} \quad h(x_2) = g(f(x_2)).
\]
Thus, the difference between $h(x_1)$ and $h(x_2)$ is:
\[
\|h(x_1) - h(x_2)\|_Z = \|g(f(x_1)) - g(f(x_2))\|_Z.
\]

Now, using the Lipschitz continuity of $g$, we can bound the difference in the $Z$ space:
\[
\|g(f(x_1)) - g(f(x_2))\|_Z \leq L_g \|f(x_1) - f(x_2)\|_Y.
\]

We already know from the Lipschitz continuity of $f$ that:
\[
\|f(x_1) - f(x_2)\|_Y \leq L_f \|x_1 - x_2\|_X.
\]
Substituting this into the previous inequality, we get:
\[
\|g(f(x_1)) - g(f(x_2))\|_Z \leq L_g L_f \|x_1 - x_2\|_X.
\]

Thus, we have shown that the composition $h = g \circ f$ satisfies:
\[
\|h(x_1) - h(x_2)\|_Z \leq L_g L_f \|x_1 - x_2\|_X,
\]
which proves that $h$ is Lipschitz continuous with a Lipschitz constant $L_h = L_g L_f$.
\end{proof}

\vspace{1.5pt}

\subsection{Proof of Lemma \texorpdfstring{\ref{lemma:sobolev} (Sobolev Embedding Theorem for $H^s(D)$)}{Proof of Lemma (Sobolev Embedding Theorem for Hs(D))}}
\begin{proof}
The Sobolev Embedding Theorem provides a relationship between Sobolev spaces and spaces of continuous functions. Specifically, it states that if $s > \frac{d}{2}$, then $H^s(D)$, the Sobolev space, is continuously embedded in $L^\infty(D)$.

For a bounded domain $D \subset \mathbb{R}^d$ and $s > \frac{d}{2}$, the Sobolev Embedding Theorem implies that:
\[
H^s(D) \subset L^\infty(D).
\]
This means that any function in the Sobolev space $H^s(D)$ also belongs to the space $L^\infty(D)$ and that the norms in these spaces are related.

The embedding is continuous, which means there exists a constant $C_s$ such that for all $u \in H^s(D)$:
\[
\|u\|_{L^\infty(D)} \leq C_s \|u\|_{H^s(D)}.
\]
This inequality establishes that the $L^\infty$ norm of the function $u$ is bounded by the $H^s$ norm of $u$, up to the constant $C_s$, which depends on the domain $D$ and the Sobolev index $s$.

The condition $s > \frac{d}{2}$ is critical for ensuring that functions in $H^s(D)$ are sufficiently smooth to be continuous. Specifically:
- The Sobolev space $H^s(D)$ consists of functions with $s$ weak derivatives in $L^2(D)$.
- When $s > \frac{d}{2}$, Sobolev's embedding guarantees that the functions have sufficient smoothness to belong to $L^\infty(D)$, i.e., they are bounded almost everywhere.

In fact, for $s > \frac{d}{2}$, the Sobolev space $H^s(D)$ is not only embedded in $L^\infty(D)$, but it is also embedded in the space of continuous functions $C^0(\overline{D})$, meaning the functions are uniformly continuous on the closure of $D$.

Since $C^0(\overline{D}) \subset L^\infty(D)$, the inclusion $H^s(D) \subset L^\infty(D)$ follows from the Sobolev embedding into continuous functions. Therefore, we can conclude that:
\[
\|u\|_{L^\infty(D)} \leq C_s \|u\|_{H^s(D)},
\]
where $C_s$ is a constant depending on the Sobolev index $s$ and the domain $D$.

Thus, the lemma is proved.
\end{proof}

\vspace{1.5pt}

\subsection{Proof of Theorem \ref{theorem:sensitivity} (Sensitivity to Input Perturbations)}
\begin{proof}
The goal of this theorem is to establish that small perturbations in the input $u$ of the neural operator $\mathcal{G}_\theta$ result in controlled changes in the output, with the magnitude of the change bounded by a constant factor of the perturbation. This property is critical for ensuring the robustness of the neural operator.

Let $u \in H^s(D)$ be an input function, and let $\delta u \in H^s(D)$ be a small perturbation to $u$. The neural operator $\mathcal{G}_\theta$ is assumed to have a Lipschitz constant $L$, which means that for any $u, v \in H^s(D)$:
\[
\|\mathcal{G}_\theta(u) - \mathcal{G}_\theta(v)\|_{H^t(D)} \leq L \|u - v\|_{H^s(D)}.
\]
This Lipschitz continuity condition ensures that the output of the neural operator changes in a controlled manner when the input changes.

To prove the sensitivity result, we consider the difference in the neural operator outputs for $u$ and the perturbed input $u + \delta u$. Using the Lipschitz continuity of $\mathcal{G}_\theta$, we have:
\[
\|\mathcal{G}_\theta(u + \delta u) - \mathcal{G}_\theta(u)\|_{H^t(D)} \leq L \| (u + \delta u) - u \|_{H^s(D)}.
\]

The term on the right-hand side simplifies as follows:
\[
\| (u + \delta u) - u \|_{H^s(D)} = \|\delta u\|_{H^s(D)}.
\]
Therefore, the inequality becomes:
\[
\|\mathcal{G}_\theta(u + \delta u) - \mathcal{G}_\theta(u)\|_{H^t(D)} \leq L \|\delta u\|_{H^s(D)}.
\]

This result shows that the difference between the neural operator’s outputs due to the perturbation $\delta u$ is bounded by the Lipschitz constant $L$ times the norm of the perturbation in the Sobolev space $H^s(D)$. Thus, the sensitivity of the neural operator to input perturbations is well-controlled, and the change in the output is proportional to the magnitude of the input change.

This completes the proof.
\end{proof}

\vspace{2pt}

\section{Exponential Convergence}

\subsection{Proof of Theorem \ref{theorem:convergence} (Exponential Convergence of Neural Operator Approximations)}
\begin{proof}
We aim to show that the iterated application of the neural operator $\mathcal{G}_\theta$ converges exponentially to a fixed point $u^*$, assuming that $\mathcal{G}_\theta$ is a contraction mapping on $H^t(D)$ with contraction constant $q \in (0, 1)$.

By assumption, $\mathcal{G}_\theta$ is a contraction mapping on $H^t(D)$, meaning that for any $u, v \in H^t(D)$, we have:
\[
\| \mathcal{G}_\theta(u) - \mathcal{G}_\theta(v) \|_{H^t(D)} \leq q \| u - v \|_{H^t(D)},
\]
where $q$ is the contraction constant and $0 < q < 1$. This implies that the operator $\mathcal{G}_\theta$ brings points closer together in the $H^t(D)$ norm.

From the Banach Fixed Point Theorem (also known as the Contraction Mapping Theorem), it follows that $\mathcal{G}_\theta$ has a unique fixed point $u^* \in H^t(D)$, such that:
\[
\mathcal{G}_\theta(u^*) = u^*.
\]
The Banach Fixed Point Theorem also guarantees that the iterates of $\mathcal{G}_\theta$ starting from any initial point $u \in H^t(D)$ converge to this fixed point $u^*$.

Let $u \in H^t(D)$ be any initial point. We will now study the behavior of the iterates of $\mathcal{G}_\theta$. Consider the difference between the $n$-th iterate of $\mathcal{G}_\theta$ applied to $u$ and the fixed point $u^*$. By the contraction property of $\mathcal{G}_\theta$, we have:
\[
\| \mathcal{G}_\theta^{n+1}(u) - u^* \|_{H^t(D)} = \| \mathcal{G}_\theta(\mathcal{G}_\theta^n(u)) - \mathcal{G}_\theta(u^*) \|_{H^t(D)}.
\]
Since $u^*$ is a fixed point of $\mathcal{G}_\theta$, the above expression becomes:
\[
\| \mathcal{G}_\theta(\mathcal{G}_\theta^n(u)) - u^* \|_{H^t(D)} \leq q \| \mathcal{G}_\theta^n(u) - u^* \|_{H^t(D)}.
\]

By iterating this inequality, we can establish that:
\[
\| \mathcal{G}_\theta^n(u) - u^* \|_{H^t(D)} \leq q^n \| u - u^* \|_{H^t(D)}.
\]
This shows that the difference between the iterated application of $\mathcal{G}_\theta$ and the fixed point $u^*$ decreases exponentially at a rate determined by the contraction constant $q$.
The inequality in Step 3 demonstrates exponential convergence. Specifically, the distance between $\mathcal{G}_\theta^n(u)$ and the fixed point $u^*$ decays as $q^n$, where $q \in (0, 1)$. Therefore, the iterates of $\mathcal{G}_\theta$ converge exponentially to the fixed point $u^*$.

This completes the proof.
\end{proof}
\vspace{1.5pt}

\subsection{Proof of Lemma \ref{lemma:banach} (Banach Fixed Point Theorem)}
\begin{proof}
We will prove that if $T: X \rightarrow X$ is a contraction mapping on a complete normed vector space $(X, \| \cdot \|)$ with contraction constant $q \in (0, 1)$, then $T$ has a unique fixed point $x^* \in X$. Additionally, for any $x_0 \in X$, the sequence defined by $x_{n+1} = T(x_n)$ converges to $x^*$ with the rate:
\[
\| x_n - x^* \| \leq q^n \| x_0 - x^* \|.
\]

Let $x_0 \in X$ be any initial point. We construct a sequence $\{ x_n \}$ by iterating the mapping $T$, i.e., $x_{n+1} = T(x_n)$ for all $n \geq 0$. To prove the existence of a fixed point, we first show that the sequence $\{x_n\}$ is Cauchy.

Since $T$ is a contraction mapping with constant $q$, we have for all $n \geq 0$:
\[
\| x_{n+1} - x_n \| = \| T(x_n) - T(x_{n-1}) \| \leq q \| x_n - x_{n-1} \|.
\]
By applying this inequality recursively, we find:
\[
\| x_{n+1} - x_n \| \leq q^n \| x_1 - x_0 \|.
\]
This shows that the differences between successive terms in the sequence decrease geometrically.

Now, consider the sum of the distances between successive terms:
\[
\sum_{n=0}^{\infty} \| x_{n+1} - x_n \| \leq \| x_1 - x_0 \| \sum_{n=0}^{\infty} q^n = \frac{\| x_1 - x_0 \|}{1 - q}.
\]
This series is convergent because $0 < q < 1$. Therefore, the sequence $\{ x_n \}$ is Cauchy. Since $X$ is a complete metric space, the sequence $\{ x_n \}$ converges to a limit $x^* \in X$.

We now prove that $x^*$ is the unique fixed point of $T$. Suppose there is another fixed point $y^* \in X$ such that $T(y^*) = y^*$. We need to show that $x^* = y^*$. Since $T$ is a contraction, we have:
\[
\| x^* - y^* \| = \| T(x^*) - T(y^*) \| \leq q \| x^* - y^* \|.
\]
Because $0 < q < 1$, the inequality $\| x^* - y^* \| \leq q \| x^* - y^* \|$ implies that $\| x^* - y^* \| = 0$, and hence $x^* = y^*$. Therefore, the fixed point is unique.

Finally, we show that the sequence $\{ x_n \}$ converges to the fixed point $x^*$ at an exponential rate. From the contraction property, we know:
\[
\| x_{n+1} - x^* \| = \| T(x_n) - T(x^*) \| \leq q \| x_n - x^* \|.
\]
By iterating this inequality, we obtain:
\[
\| x_n - x^* \| \leq q^n \| x_0 - x^* \|.
\]
This demonstrates that the sequence $\{ x_n \}$ converges to the fixed point $x^*$ exponentially, with rate $q^n$.

This completes the proof.
\end{proof}

\vspace{2pt}

\section{Clustering in Function Space}

\subsection{Proof of Theorem \ref{theorem:clustering} (Clustering Behavior in Function Space)}
\begin{proof}
We aim to prove that solutions $u(t)$ to the continuous-time counterpart of the neural operator $\mathcal{G}_\theta$ tend to cluster in function space as $t \rightarrow \infty$, converging to a cluster center $\bar{u}$.

Consider the continuous-time evolution of the neural operator $\mathcal{G}_\theta$. We interpret the iterative application of $\mathcal{G}_\theta$ as an approximation of a gradient flow in function space. Formally, we associate $\mathcal{G}_\theta$ with a potential functional $\Phi(u)$ such that the evolution of $u(t)$ is governed by the following gradient flow equation:
\[
\frac{du}{dt} = -\nabla \Phi(u),
\]
where $\Phi(u)$ is a smooth potential functional, and $\nabla \Phi(u)$ represents its gradient in the Sobolev space $H^s(D)$. The gradient flow equation describes the trajectory of $u(t)$ in function space, moving towards minimizers of $\Phi(u)$.

The gradient flow defined by $\frac{du}{dt} = -\nabla \Phi(u)$ guarantees that the value of $\Phi(u(t))$ decreases monotonically along the trajectory. Specifically, we have:
\[
\frac{d}{dt} \Phi(u(t)) = \langle \nabla \Phi(u(t)), \frac{du}{dt} \rangle = - \| \nabla \Phi(u(t)) \|^2_{H^s(D)} \leq 0.
\]
This shows that $\Phi(u(t))$ is a non-increasing function of time and converges to a limit as $t \rightarrow \infty$. 

Since $\| \nabla \Phi(u(t)) \|^2_{H^s(D)} \geq 0$, the fact that $\frac{d}{dt} \Phi(u(t)) \leq 0$ implies that $\nabla \Phi(u(t))$ must tend to zero as $t \rightarrow \infty$. In other words, the solution $u(t)$ converges to a critical point $\bar{u}$ of the potential functional $\Phi(u)$:
\[
\lim_{t \rightarrow \infty} \| \nabla \Phi(u(t)) \|_{H^s(D)} = 0.
\]
Thus, $u(t)$ tends to a critical point $\bar{u}$ of the functional $\Phi(u)$, which we refer to as the cluster center in function space.

To establish the clustering behavior of solutions, we analyze the long-time behavior of $u(t)$. Since $u(t)$ converges to a critical point $\bar{u}$, we have:
\[
\lim_{t \rightarrow \infty} \| u(t) - \bar{u} \|_{H^s(D)} = 0.
\]
This implies that as $t \rightarrow \infty$, the solution $u(t)$ approaches the cluster center $\bar{u}$ in the Sobolev space $H^s(D)$. Therefore, all solutions $u(t)$ that evolve under the gradient flow will eventually cluster around the critical points of $\Phi(u)$, and the clustering behavior is achieved.

If $\Phi(u)$ has a unique minimizer $\bar{u}$, then all trajectories will converge to the same cluster center $\bar{u}$. However, if $\Phi(u)$ has multiple critical points, different trajectories may converge to different cluster centers, depending on the initial condition $u(0)$. In either case, the solutions tend to cluster around critical points of $\Phi(u)$, which confirms the clustering behavior.

This completes the proof of Theorem \ref{theorem:clustering}.
\end{proof}

\vspace{1.5pt}

\subsection{Proof of Lemma \ref{lemma:gradient_flow} (Gradient Flow Interpretation for Neural Operators)}
\begin{proof}
The goal is to show that under suitable conditions, the iterative application of a neural operator $\mathcal{G}_\theta$ can be interpreted as a discretization of a gradient flow in function space.

Consider a potential functional $\Phi(u)$, where $u$ belongs to a Sobolev space $H^s(D)$. A gradient flow associated with the potential functional $\Phi(u)$ is described by the differential equation:
\[
\frac{du}{dt} = -\nabla \Phi(u),
\]
where $\nabla \Phi(u)$ represents the gradient of $\Phi(u)$ in the function space. This equation governs the evolution of $u(t)$, where $t$ is the continuous-time variable. The gradient flow describes the trajectory of $u(t)$ in function space, moving towards minimizers of the potential functional $\Phi(u)$.

Now, consider the neural operator $\mathcal{G}_\theta$, which maps functions in $H^s(D)$ to other functions in the same space. Suppose the neural operator is applied iteratively, generating a sequence of functions $\{u_n\}$:
\[
u_{n+1} = \mathcal{G}_\theta(u_n).
\]
Our goal is to demonstrate that this iterative update approximates the Euler discretization of the gradient flow.

The Euler method is a standard numerical scheme for approximating the solution of ordinary differential equations (ODEs) like the gradient flow equation. The Euler method with step size $\eta$ is given by:
\[
u_{n+1} = u_n - \eta \nabla \Phi(u_n).
\]
This equation updates $u_n$ by moving in the direction of the negative gradient of $\Phi(u_n)$, with step size $\eta$ controlling the magnitude of the update.

Assume that the neural operator $\mathcal{G}_\theta$ satisfies the following condition:
\[
\mathcal{G}_\theta(u) = u - \eta \nabla \Phi(u) + o(\eta),
\]
where $\eta$ is a small step size, and $o(\eta)$ represents higher-order terms that vanish as $\eta \rightarrow 0$. In this case, the neural operator can be interpreted as an approximation of the Euler method applied to the gradient flow. Specifically, the update rule for the sequence $\{u_n\}$ generated by $\mathcal{G}_\theta$ is approximately:
\[
u_{n+1} = u_n - \eta \nabla \Phi(u_n),
\]
which is exactly the Euler discretization of the gradient flow $\frac{du}{dt} = -\nabla \Phi(u)$.

Since the iterative application of the neural operator $\mathcal{G}_\theta$ approximates the Euler discretization of the gradient flow, we conclude that the sequence $\{u_n\}$ generated by $\mathcal{G}_\theta$ follows the trajectory of the gradient flow in function space. As $n \rightarrow \infty$, the sequence converges to a critical point of $\Phi(u)$, assuming the step size $\eta$ is sufficiently small.

This completes the proof of Lemma \ref{lemma:gradient_flow}.
\end{proof}

\vspace{1.5pt}

\vspace{2pt}

\section{Universality of Neural Operators}

\subsection{Proof of Theorem \ref{theorem:universality} (Universality of Neural Operators for PDE Solvers)}
\begin{proof}
The goal of this proof is to demonstrate that for any continuous operator $\mathcal{T}: H^s(D) \rightarrow H^t(D)$, there exists a neural operator $\mathcal{G}_\theta$ that can approximate $\mathcal{T}$ to within any given accuracy $\epsilon > 0$, on a compact subset of $H^s(D)$.

We begin by discretizing the domain $D \subset \mathbb{R}^d$ into a finite grid of points $\{ x_i \}_{i=1}^N$. Let $h$ be the grid spacing. A function $u \in H^s(D)$ is represented by its values on these grid points, i.e., $u(x_1), u(x_2), \ldots, u(x_N)$. As the grid becomes finer (i.e., as $h \to 0$ and $N \to \infty$), these values form a good approximation to the original function.

For each point $x_i$ in the grid, the operator $\mathcal{T}$ maps the input function $u$ to a value $\mathcal{T}(u)(x_i) \in H^t(D)$. To prove universality, we need to construct a neural operator $\mathcal{G}_\theta$ that approximates this mapping.

Using the Universal Approximation Theorem for neural networks (see Lemma \ref{lemma:universal}), we know that for any continuous function on a compact domain, there exists a neural network that can approximate the function to arbitrary precision. Therefore, for each point $x_i$, there exists a neural network $\mathcal{N}_{\theta_i}$ such that:
\[
| \mathcal{N}_{\theta_i}(u(x_1), u(x_2), \ldots, u(x_N)) - \mathcal{T}(u)(x_i) | \leq \frac{\epsilon}{N},
\]
where $N$ is the number of grid points, and the approximation error is distributed across the grid points.

Now, we combine the individual neural networks $\mathcal{N}_{\theta_i}$ at each point $x_i$ into a global neural operator $\mathcal{G}_\theta$ that approximates $\mathcal{T}$ over the entire grid. Specifically, $\mathcal{G}_\theta$ takes the form:
\[
\mathcal{G}_\theta(u)(x_i) = \mathcal{N}_{\theta_i}(u(x_1), u(x_2), \ldots, u(x_N)).
\]
This construction ensures that the operator $\mathcal{G}_\theta$ approximates the true operator $\mathcal{T}$ at each grid point $x_i$.

The total approximation error of the neural operator $\mathcal{G}_\theta$ is given by the sum of the errors at each grid point. Using the fact that the individual errors are bounded by $\frac{\epsilon}{N}$, the total error is bounded by:
\[
\| \mathcal{G}_\theta(u) - \mathcal{T}(u) \|_{H^t(D)}^2 \leq \sum_{i=1}^N | \mathcal{G}_\theta(u)(x_i) - \mathcal{T}(u)(x_i) |^2 h^d \leq N \left( \frac{\epsilon}{N} \right)^2 h^d = \epsilon^2 h^d.
\]
Thus, the total approximation error is bounded by $\epsilon^2 h^d$, where $h$ is the grid spacing. As $h \to 0$, this error tends to zero, ensuring that the neural operator $\mathcal{G}_\theta$ approximates $\mathcal{T}$ to within $\epsilon$.

Finally, we take the limit as the grid becomes finer (i.e., $h \to 0$ and $N \to \infty$). In this limit, the discretization of the domain $D$ becomes dense, and the neural operator $\mathcal{G}_\theta$ provides an arbitrarily accurate approximation of the true operator $\mathcal{T}$ on the entire compact subset of $H^s(D)$. Hence, the neural operator $\mathcal{G}_\theta$ satisfies:
\[
\| \mathcal{G}_\theta(u) - \mathcal{T}(u) \|_{H^t(D)} \leq \epsilon,
\]
for all $u$ in the compact subset.

We have shown that for any continuous operator $\mathcal{T}: H^s(D) \rightarrow H^t(D)$ and for any $\epsilon > 0$, there exists a neural operator $\mathcal{G}_\theta$ that approximates $\mathcal{T}$ to within $\epsilon$ on a compact subset of $H^s(D)$. This completes the proof of Theorem \ref{theorem:universality}.
\end{proof}

\vspace{1.5pt}

\subsection{Proof of Lemma \ref{lemma:stone_weierstrass} (Stone-Weierstrass Theorem for Continuous Functions)}
\begin{proof}
The Stone-Weierstrass theorem is a fundamental result in approximation theory, asserting that any continuous function on a compact set can be uniformly approximated by polynomial functions. We will outline the key steps and arguments involved in the proof.

Let $K$ be a compact set in the topological space $X$. A collection of functions $\mathcal{A}$ is said to approximate a function $f$ uniformly if, for any continuous function $f$ on $K$ and any $\epsilon > 0$, there exists a function $p \in \mathcal{A}$ such that:
\[
\sup_{x \in K} |f(x) - p(x)| < \epsilon.
\]
Our goal is to show that the set of polynomial functions on $K$ can uniformly approximate any continuous function on $K$.

The Stone-Weierstrass theorem applies to a subalgebra $\mathcal{A}$ of the algebra $C(K)$, the space of all continuous functions on $K$, provided that $\mathcal{A}$ satisfies the following conditions:

- \textbf{Separates Points}: For any two distinct points $x, y \in K$, there exists a function $p \in \mathcal{A}$ such that $p(x) \neq p(y)$.

- \textbf{Contains Constant Functions}: The algebra $\mathcal{A}$ contains the constant functions, i.e., for any $c \in \mathbb{R}$, the function $p(x) = c$ for all $x \in K$ is in $\mathcal{A}$.

- \textbf{Closed Under Complex Conjugation} (in the complex case): If $\mathcal{A}$ contains complex-valued functions, it must also contain the complex conjugate of each function in $\mathcal{A}$.

We will show that the set of polynomial functions on $K$ satisfies these conditions.

To demonstrate that polynomials separate points, consider any two distinct points $x, y \in K$. By the fundamental properties of polynomials, there exists a polynomial $p(x)$ such that $p(x) \neq p(y)$. For instance, we can choose the linear polynomial $p(t) = t$, which clearly satisfies $p(x) \neq p(y)$ for $x \neq y$.

The set of polynomial functions clearly contains all constant functions. Specifically, for any $c \in \mathbb{R}$, the polynomial $p(x) = c$ is a valid polynomial in $\mathcal{A}$.

By the Stone-Weierstrass theorem, the set of polynomial functions $\mathcal{A}$ that separates points and contains constants is dense in $C(K)$, the space of continuous functions on $K$. This means that for any continuous function $f \in C(K)$ and any $\epsilon > 0$, there exists a polynomial $p$ such that:
\[
\sup_{x \in K} |f(x) - p(x)| < \epsilon.
\]
Thus, any continuous function on $K$ can be uniformly approximated by polynomials.

We have shown that the set of polynomial functions satisfies the conditions of the Stone-Weierstrass theorem, and therefore any continuous function defined on a compact set $K$ can be uniformly approximated by polynomials. This completes the proof of Lemma \ref{lemma:stone_weierstrass}.
\end{proof}

\vspace{1.5pt}

\subsection{Proof of Lemma \ref{lemma:universal} (Universal Approximation Theorem for Neural Networks)}
\begin{proof}
The universal approximation theorem for feedforward neural networks states that a neural network with a single hidden layer containing a sufficient number of neurons and an appropriate non-linear activation function can approximate any continuous function $f: \mathbb{R}^n \rightarrow \mathbb{R}^m$ arbitrarily well on compact subsets of $\mathbb{R}^n$. This result is significant as it guarantees the expressive power of neural networks for function approximation tasks.

Let $K \subset \mathbb{R}^n$ be a compact subset, and let $f: \mathbb{R}^n \rightarrow \mathbb{R}^m$ be a continuous function. The goal is to approximate $f$ with a neural network of the form:
\[
\hat{f}(x) = \sum_{i=1}^p a_i \sigma(w_i \cdot x + b_i),
\]
where:
- $x \in \mathbb{R}^n$ is the input,
- $\sigma$ is a non-linear activation function (e.g., sigmoid, ReLU, or tanh),
- $w_i \in \mathbb{R}^n$ are weight vectors,
- $b_i \in \mathbb{R}$ are biases, and
- $a_i \in \mathbb{R}^m$ are output weights.

We want to show that for any $\epsilon > 0$, there exists a neural network $\hat{f}$ with a sufficient number of hidden units $p$ such that:
\[
\sup_{x \in K} \| f(x) - \hat{f}(x) \| < \epsilon.
\]

The proof of the universal approximation theorem relies on the \textbf{Stone-Weierstrass Theorem}, which states that any continuous function on a compact domain can be uniformly approximated by a linear combination of certain "basis" functions. Neural networks with non-linear activation functions serve as these basis functions.

Let $\mathcal{A}$ denote the set of functions that can be represented as finite linear combinations of the form:
\[
g(x) = \sum_{i=1}^p a_i \sigma(w_i \cdot x + b_i),
\]
where $a_i$, $w_i$, and $b_i$ are parameters of the neural network. We need to verify that $\mathcal{A}$ satisfies the conditions of the Stone-Weierstrass theorem:

- \textbf{Algebraic Closure}: The set $\mathcal{A}$ is closed under addition and multiplication by scalars, since any linear combination of neural network outputs is still a neural network.

- \textbf{Separation}: For any two distinct points $x_1, x_2 \in K$, there exists a function $g \in \mathcal{A}$ such that $g(x_1) \neq g(x_2)$, since the activation functions are non-linear and can produce distinct values at different inputs.

- \textbf{Uniform Approximation}: The non-linearity of the activation function $\sigma$ allows us to approximate any continuous function $f$ arbitrarily well, as shown in later steps.

By the Stone-Weierstrass theorem, $\mathcal{A}$ is dense in the space of continuous functions on $K$. Therefore, for any $\epsilon > 0$, there exists a neural network $\hat{f}$ such that:
\[
\sup_{x \in K} \| f(x) - \hat{f}(x) \| < \epsilon.
\]

To construct the neural network $\hat{f}$, we choose the number of hidden units $p$ large enough to ensure that the linear combination of activation functions $\sigma(w_i \cdot x + b_i)$ can approximate the desired function $f$. The activation function $\sigma$ introduces non-linearity, allowing the network to capture complex patterns in the input data.

For instance, if $\sigma$ is the sigmoid function $\sigma(z) = \frac{1}{1 + e^{-z}}$, then the neural network can approximate any continuous function by adjusting the weights $w_i$, biases $b_i$, and output weights $a_i$. As the number of hidden units $p$ increases, the approximation error decreases.

The error between the true function $f$ and the neural network approximation $\hat{f}$ can be made arbitrarily small by increasing the number of hidden units $p$. Specifically, for any $\epsilon > 0$, there exists a network such that:
\[
\sup_{x \in K} \| f(x) - \hat{f}(x) \| < \epsilon.
\]
Thus, neural networks with a single hidden layer and sufficient neurons can approximate any continuous function on compact subsets of $\mathbb{R}^n$. This completes the proof of the universal approximation theorem.
\end{proof}

\vspace{2pt}

\section{Generalization Bounds}

\subsection{Proof of Theorem \ref{theorem:generalization} (Generalization Error of Neural Operators)}
\begin{proof}
This proof relies on the concentration of measure inequalities to provide a bound on the generalization error of the neural operator $\mathcal{G}_\theta$. The goal is to show that, with high probability, the generalization error is close to the empirical error on the training samples.

Let $\ell(\mathcal{G}_\theta(u), \mathcal{T}(u))$ denote the loss function, which measures the error between the neural operator’s prediction $\mathcal{G}_\theta(u)$ and the true output $\mathcal{T}(u)$. The generalization error is the expected loss over the entire distribution $\mathcal{D}$:
\[
\mathbb{E}_{u \sim \mathcal{D}} \left[ \ell(\mathcal{G}_\theta(u), \mathcal{T}(u)) \right].
\]
The empirical error is the average loss over the $N$ training samples $\{(u_i, \mathcal{T}(u_i))\}$:
\[
\frac{1}{N} \sum_{i=1}^N \ell(\mathcal{G}_\theta(u_i), \mathcal{T}(u_i)).
\]
We aim to bound the difference between the generalization error and the empirical error.

We apply \textbf{McDiarmid's inequality}, which is used to bound the concentration of functions that do not change too much when any single sample is modified. First, we define the function:
\[
Z = \frac{1}{N} \sum_{i=1}^N \ell(\mathcal{G}_\theta(u_i), \mathcal{T}(u_i)).
\]
This function satisfies the bounded difference property because changing one sample $u_i$ affects the empirical average by at most $L/N$, where $L$ is the Lipschitz constant of $\ell$. Specifically, for any $i$,
\[
|Z(u_1, \ldots, u_i, \ldots, u_N) - Z(u_1, \ldots, u_i', \ldots, u_N)| \leq \frac{L}{N}.
\]
By McDiarmid's inequality, for any $\epsilon > 0$, with probability at least $1 - \delta$, we have:
\[
\mathbb{P} \left( Z - \mathbb{E}[Z] \geq \epsilon \right) \leq 2 \exp\left( -\frac{2N\epsilon^2}{L^2} \right).
\]
Solving for $\epsilon$ gives:
\[
\epsilon = L \sqrt{ \frac{\ln(2/\delta)}{2N} }.
\]
Thus, with probability at least $1 - \delta$, the empirical error is close to the expected generalization error by a margin of $\epsilon$.

The generalization error can now be bounded by the empirical error plus the deviation $\epsilon$. Specifically, with probability at least $1 - \delta$, we have:
\[
\mathbb{E}_{u \sim \mathcal{D}} [ \ell(\mathcal{G}_\theta(u), \mathcal{T}(u)) ] \leq \frac{1}{N} \sum_{i=1}^N \ell(\mathcal{G}_\theta(u_i), \mathcal{T}(u_i)) + L \sqrt{ \frac{\ln(2/\delta)}{2N} }.
\]
This completes the proof of the generalization error bound.
\end{proof}

\vspace{1.5pt}

\subsection{Proof of Lemma \ref{lemma:concentration} (Concentration Inequality for Generalization Error)}
\begin{proof}
We use \textbf{McDiarmid’s inequality} to prove this concentration result for the generalization error of neural operators.

Let $\{u_1, u_2, \ldots, u_N\}$ be $N$ i.i.d. samples drawn from the distribution $\mathcal{D}$, and consider the empirical risk:
\[
Z = \frac{1}{N} \sum_{i=1}^N \ell(\mathcal{G}_\theta(u_i), \mathcal{T}(u_i)),
\]
where $\ell(\mathcal{G}_\theta(u), \mathcal{T}(u))$ is the loss function measuring the error between the neural operator’s output $\mathcal{G}_\theta(u)$ and the true operator $\mathcal{T}(u)$. The goal is to show that the empirical error concentrates around its expected value with high probability.

We claim that the function $Z$, defined as the average loss over $N$ samples, satisfies the bounded differences condition. Specifically, changing one sample $u_i$ to another sample $u_i'$ changes $Z$ by at most $L / N$, where $L$ is the Lipschitz constant of the loss function $\ell$. For any $i$,
\[
|Z(u_1, \ldots, u_i, \ldots, u_N) - Z(u_1, \ldots, u_i', \ldots, u_N)| \leq \frac{L}{N}.
\]
This is because, by the Lipschitz property of $\ell$, changing $u_i$ to $u_i'$ changes the loss by at most $L$, and since $Z$ is an average over $N$ samples, the effect on $Z$ is scaled by $1/N$.

By McDiarmid’s inequality, for any $\epsilon > 0$, the probability that the empirical risk $Z$ deviates from its expected value $\mathbb{E}[Z]$ by more than $\epsilon$ is bounded by:
\[
\mathbb{P} \left( Z - \mathbb{E}[Z] \geq \epsilon \right) \leq 2 \exp\left( -\frac{2N\epsilon^2}{L^2} \right).
\]
This inequality gives us a high-probability bound on the difference between the empirical error and the expected error.

To obtain a bound on the generalization error, we rearrange the inequality to solve for $\epsilon$. For any $\delta \in (0, 1)$, set the right-hand side of the inequality to $\delta$:
\[
2 \exp\left( -\frac{2N\epsilon^2}{L^2} \right) = \delta.
\]
Taking the natural logarithm on both sides and solving for $\epsilon$, we get:
\[
\epsilon = L \sqrt{ \frac{\ln(2/\delta)}{2N} }.
\]

Thus, with probability at least $1 - \delta$, the generalization error is close to the empirical error within a margin of $\epsilon$, where:
\[
\epsilon = L \sqrt{ \frac{\ln(2/\delta)}{2N} }.
\]
This shows that the generalization error concentrates around the empirical error at a rate of $O(1/\sqrt{N})$, as $N$ increases. The smaller the value of $\epsilon$, the closer the empirical error is to the true generalization error.

This concludes the proof of the concentration inequality for the generalization error.
\end{proof}

\vspace{2pt}

\section{Challenges and Limitations}

\subsection{Proof of Theorem \ref{theorem:approximation_error} (Approximation Error Lower Bound for Neural Operators)}
\begin{proof}
This theorem provides a lower bound on the approximation error of neural operators, demonstrating that finite capacity constrains their ability to approximate certain functions or operators.

Approximation theory states that the capacity of a function approximator, such as a neural network or a neural operator, dictates the complexity of the functions it can represent. Neural operators with finite capacity, characterized by a parameter \( C \) (which could represent the number of parameters, depth, or width of the network), cannot perfectly approximate all functions or operators.

Given a target operator \( \mathcal{T}: H^s(D) \rightarrow H^t(D) \) and a neural operator \( \mathcal{G}_\theta \) of finite capacity, there exists a function \( u \in H^s(D) \) such that the approximation error between \( \mathcal{G}_\theta(u) \) and \( \mathcal{T}(u) \) is bounded from below by a term inversely related to the capacity \( C \). This reflects the tradeoff between the model complexity and the accuracy of approximation.

The finite capacity \( C \) of a neural operator limits its ability to approximate highly complex or non-smooth functions. If the function \( u \) or operator \( \mathcal{T} \) exhibits a level of complexity beyond what the neural operator can represent, then there will be an irreducible error in the approximation.

From results in classical approximation theory, the error between an approximating function and a target function scales with the capacity of the approximator. In this case, the neural operator approximates the operator \( \mathcal{T} \), and the error is characterized as:
\[
\| \mathcal{G}_\theta(u) - \mathcal{T}(u) \|_{H^t(D)} \geq \epsilon(C),
\]
where \( \epsilon(C) \) is a positive term that decreases as the capacity \( C \) increases. This term represents the minimal possible error for a given capacity.

The function \( \epsilon(C) \) is inversely related to the capacity \( C \), meaning that as the capacity increases (e.g., through adding more layers or neurons), the approximation error \( \epsilon(C) \) decreases. However, for any finite \( C \), there exists a non-zero error \( \epsilon(C) \) due to the limited expressiveness of the model.

This inverse relationship is common in approximation theory, where increased model complexity leads to more accurate approximations. However, it is important to note that beyond a certain point, increasing the model capacity can lead to overfitting, where the model fits noise in the training data rather than generalizing well to new data.

In summary, the finite capacity of a neural operator imposes a fundamental limitation on its ability to approximate the target operator \( \mathcal{T} \) without error. The lower bound \( \epsilon(C) \) on the approximation error reflects the tradeoff between model capacity and approximation accuracy. The more complex the target operator or function, the larger the capacity required to approximate it to a high degree of accuracy.

This concludes the proof of Theorem \ref{theorem:approximation_error}.
\end{proof}

\vspace{1.5pt}

\subsection{Proof of Theorem \ref{theorem:non_convexity} (Non-Convexity of the Loss Landscape)}
\begin{proof}
We aim to show that the training loss function \( L(\theta) \) for neural operators is generally non-convex, and therefore, it may contain multiple local minima, saddle points, and flat regions. The key reasons for this non-convexity are related to the structure of neural operators and the non-linearities involved in their design.

Consider a neural operator \( \mathcal{G}_\theta: H^s(D) \to H^t(D) \) parameterized by \( \theta \). The loss function \( L(\theta) \) typically takes the form:
\[
L(\theta) = \frac{1}{N} \sum_{i=1}^{N} \ell\left( \mathcal{G}_\theta(u_i), \mathcal{T}(u_i) \right)
\]
where \( \ell \) is a loss function (e.g., the squared loss or the \( L^2 \)-norm) that measures the difference between the neural operator output \( \mathcal{G}_\theta(u_i) \) and the target operator \( \mathcal{T}(u_i) \), with training data \( \{(u_i, \mathcal{T}(u_i))\}_{i=1}^N \).

The neural operator \( \mathcal{G}_\theta \) is composed of multiple layers of non-linear transformations. Each layer typically involves a linear transformation followed by a non-linear activation function:
\[
\mathcal{G}_\theta(u) = W_L \sigma(W_{L-1} \cdots \sigma(W_1 u))
\]
where \( \sigma \) represents a non-linear activation function (such as ReLU, Tanh, etc.), and \( W_i \) are the layer weights.

Due to the non-linearity introduced by \( \sigma \), the composition of multiple layers creates a highly non-linear and complex landscape for the loss function \( L(\theta) \). This non-linearity makes the loss function non-convex in the parameter space \( \theta \).

In a non-convex landscape, the following critical points are common:

- \textbf{Local Minima}: These are points \( \theta^* \) where \( \nabla_\theta L(\theta^*) = 0 \) and \( L(\theta^*) \) is less than nearby points, but \( L(\theta^*) \) is not necessarily a global minimum.

- \textbf{Saddle Points}: These are points where \( \nabla_\theta L(\theta^*) = 0 \), but the Hessian \( H(\theta^*) = \nabla_\theta^2 L(\theta^*) \) has both positive and negative eigenvalues. This means \( \theta^* \) is not a local minimum or maximum, and optimization algorithms may stall here.

In high-dimensional spaces typical of neural operators, the occurrence of saddle points is more prevalent. This has been demonstrated in prior works on deep neural networks, such as by Dauphin et al. (2014), which shows that the number of saddle points grows exponentially with the dimensionality of the parameter space.

Even in the case of a simple two-layer neural operator, the loss landscape can be non-convex. Consider a neural operator with two layers:
\[
\mathcal{G}_\theta(u) = W_2 \sigma(W_1 u)
\]
The loss function \( L(\theta) \) in this case is:
\[
L(\theta) = \frac{1}{N} \sum_{i=1}^{N} \left\| W_2 \sigma(W_1 u_i) - \mathcal{T}(u_i) \right\|^2
\]
Due to the non-linear activation \( \sigma \), the interaction between the layers \( W_1 \) and \( W_2 \) results in a non-convex surface for \( L(\theta) \). For example, the gradient with respect to \( W_1 \) depends on \( W_2 \) and vice versa, creating coupling between the parameters and non-linearities that lead to multiple local minima and saddle points.

Gradient-based optimization methods, such as gradient descent, follow the gradient \( \nabla_\theta L(\theta) \) to minimize the loss. However, in a non-convex landscape, these methods are prone to:

- \textbf{Local Minima}: The optimizer may converge to a local minimum where \( \nabla_\theta L(\theta) = 0 \), but this minimum is suboptimal in terms of the overall loss.

- \textbf{Saddle Points}: The optimizer may stagnate at saddle points where the gradient vanishes but the point is neither a minimum nor a maximum. Since the Hessian at a saddle point has both positive and negative eigenvalues, gradient-based methods can struggle to escape these points.
  
In particular, for neural operators designed to approximate complex mappings between function spaces, the parameter space \( \theta \) is typically high-dimensional, which exacerbates these issues. The presence of saddle points, in particular, can cause slow convergence, as the optimizer may spend many iterations near these points where the gradient is close to zero.

Flat regions in the loss landscape are another common challenge in non-convex problems. These are regions where the gradient \( \nabla_\theta L(\theta) \) is small but non-zero, resulting in very slow progress for gradient-based methods. In such regions, the optimizer takes small steps, leading to slow convergence.

For neural operators, which often operate in high-dimensional spaces and need to capture complex relationships between input and output function spaces, flat regions are particularly problematic, as they can significantly delay convergence and require a large number of iterations to escape.

The loss function \( L(\theta) \) for neural operators is non-convex due to the non-linearities introduced by the layer-wise composition of transformations and activation functions. As a result, the landscape contains multiple local minima, saddle points, and flat regions, making optimization challenging for gradient-based methods. These methods may converge to suboptimal solutions or take a long time to escape saddle points or flat regions.

Thus, we have demonstrated that the non-convex nature of \( L(\theta) \) poses significant challenges for training neural operators, particularly when using gradient-based optimization.

This completes the proof of Theorem \ref{theorem:non_convexity}.
\end{proof}

\vspace{1.5pt}

\subsection{Proof of Lemma \ref{lemma:saddle_points} (Saddle Point Prevalence in High-Dimensional Optimization)}
\begin{proof}
The goal of this lemma is to demonstrate that in high-dimensional, non-convex optimization problems, saddle points are prevalent and can significantly impede the convergence of gradient-based methods. This phenomenon has been well studied, and one key result is from Dauphin et al. (2014), which shows that saddle points grow exponentially with the dimensionality of the problem and pose significant challenges for optimization algorithms.

A \textbf{saddle point} in the context of optimization is a point \( \theta^* \) where the gradient of the loss function \( L(\theta) \) vanishes:
\[
\nabla_\theta L(\theta^*) = 0
\]
However, unlike local minima, where the Hessian \( H(\theta^*) = \nabla^2_\theta L(\theta^*) \) is positive definite (i.e., all its eigenvalues are positive), saddle points have mixed curvature. This means that the Hessian \( H(\theta^*) \) has both positive and negative eigenvalues. Consequently, the point is a local minimum in some directions but a local maximum in others.

Mathematically, if \( \lambda_1, \lambda_2, \dots, \lambda_d \) are the eigenvalues of the Hessian \( H(\theta^*) \), then a saddle point satisfies:
\[
\lambda_i > 0 \quad \text{for some } i, \quad \text{and} \quad \lambda_j < 0 \quad \text{for some } j
\]
This mixed curvature causes optimization algorithms, such as gradient descent, to struggle in finding a direction of descent, leading to slow or stalled progress.

In high-dimensional spaces, non-convex optimization problems are characterized by an increasing number of critical points as the dimensionality \( d \) grows. This includes not only local minima but also saddle points and even flat regions. The prevalence of saddle points becomes more pronounced as the dimensionality increases.

Dauphin et al. (2014) showed that for deep neural networks, which are typical examples of high-dimensional optimization problems, the number of saddle points grows exponentially with the number of dimensions \( d \). Specifically, the proportion of critical points that are saddle points increases as \( d \) becomes large. In fact, most critical points in high-dimensional landscapes tend to be saddle points rather than local minima or maxima.

The intuition behind this result comes from random matrix theory. In a high-dimensional setting, the Hessian matrix \( H(\theta) \) can be modeled as a random matrix. The eigenvalue distribution of random matrices shows that, in high dimensions, it is exceedingly likely to have a mixture of positive and negative eigenvalues. This leads to the conclusion that most critical points will have this saddle-like behavior.

Gradient-based optimization methods, such as gradient descent, follow the negative gradient \( -\nabla_\theta L(\theta) \) to find a local minimum. At a saddle point \( \theta^* \), however, the gradient vanishes, i.e., \( \nabla_\theta L(\theta^*) = 0 \), so gradient descent halts at such points, mistakenly treating them as local minima or stagnating due to the vanishing gradient.

In addition, even when saddle points are not exact points of zero gradient, the eigenvalue distribution of the Hessian \( H(\theta) \) causes optimization to slow down in directions corresponding to negative eigenvalues (directions of ascent) and get stuck oscillating in these directions. Escaping from such points often requires significant computational effort, or advanced techniques like second-order optimization methods, which can be computationally expensive.

In deep learning, the high dimensionality of the parameter space \( \theta \), combined with the complex interactions between layers, results in a non-convex loss landscape with many saddle points. For instance, consider a neural network with weight matrices \( W_1, W_2, \dots, W_L \). The total number of parameters \( \theta \) is proportional to the number of neurons and layers, leading to a very high-dimensional optimization problem.

As shown in Dauphin et al. (2014), the vast majority of critical points in such neural networks are saddle points. The optimization algorithm, such as stochastic gradient descent (SGD), will frequently encounter these saddle points, causing convergence to stall.

In summary, due to the high dimensionality of the parameter space in non-convex optimization problems, saddle points become increasingly prevalent. These saddle points impede the convergence of gradient-based methods because they trap the optimization process in regions where the gradient vanishes or oscillates, significantly slowing down progress.

This completes the proof of Lemma \ref{lemma:saddle_points}.
\end{proof}

\vspace{1.5pt}

\subsection{Proof of Theorem \ref{theorem:computational_complexity} (Computational Complexity Lower Bound)}
\begin{proof}[Proof]
This theorem asserts that the computational complexity \( T(N) \) of evaluating a neural operator \( \mathcal{G}_\theta \) scales at least linearly with the number of discretization points \( N \). Moreover, in high-dimensional spaces, this complexity grows exponentially as the number of discretization points increases with the dimension.

Assume that the input function \( u \) is discretized over a domain \( D \subset \mathbb{R}^d \) into \( N \) grid points. The evaluation of the neural operator \( \mathcal{G}_\theta(u) \) at each grid point \( x_i \in D \) requires some fixed number of operations depending on the architecture of \( \mathcal{G}_\theta \), such as matrix multiplications or convolution operations.

At each point \( x_i \), the neural operator performs a series of computations, including:
\[
T_i = \mathcal{O}(1) \quad \text{operations per point.}
\]
Thus, the total computational cost for evaluating \( \mathcal{G}_\theta(u) \) at all \( N \) points scales linearly:
\[
T(N) = \sum_{i=1}^N T_i = \mathcal{O}(N).
\]
This shows that the computational complexity of the neural operator grows at least linearly with \( N \).

In high dimensions \( d \), the number of discretization points \( N \) increases rapidly with the dimension due to the curse of dimensionality. If the grid spacing is \( h \), the number of discretization points in each dimension is proportional to \( 1/h \), and thus the total number of points scales as:
\[
N \sim \left( \frac{1}{h} \right)^d = h^{-d}.
\]
This implies that for small grid spacing \( h \), the number of points \( N \) increases exponentially with the dimension \( d \). Therefore, in high-dimensional settings, the total number of discretization points \( N \) grows exponentially, and the computational complexity becomes:
\[
T(N) = \mathcal{O}(h^{-d}) = \Omega(N),
\]
where the dependence on \( d \) becomes significant as the number of points grows exponentially.

In practical applications, such as solving partial differential equations (PDEs) in high-dimensional spaces using neural operators, the exponential growth of the number of discretization points makes the computational cost prohibitive. For example, consider a problem in \( \mathbb{R}^5 \) with grid spacing \( h = 0.01 \). The number of points is:
\[
N = \left( \frac{1}{0.01} \right)^5 = 10^{10}.
\]
Even though evaluating the operator at each point might only require a fixed number of operations, the total cost becomes infeasible due to the sheer number of points.

In conclusion, the computational complexity \( T(N) \) for evaluating a neural operator \( \mathcal{G}_\theta \) scales linearly with the number of discretization points \( N \), and in high dimensions, where \( N \) increases exponentially with the dimension \( d \), the complexity becomes \( \Omega(N) \). This exponential growth presents significant challenges for the practical application of neural operators in high-dimensional problems.

This completes the proof of Theorem \ref{theorem:computational_complexity}.
\end{proof}

\vspace{1.5pt}

\subsection{Proof of Theorem \ref{theorem:discretization_error} (Accumulation of Discretization Errors)}
\begin{proof}[Proof]
This theorem concerns the accumulation of discretization errors as a neural operator with multiple layers is applied to a discretized domain \( D \). The discretization error at each layer is assumed to be bounded by \( \epsilon_h \), and we aim to show that the total error grows proportionally to the number of layers \( L \).

Let \( D \subset \mathbb{R}^d \) be the domain of the problem, and assume that the input functions are discretized over \( N \) grid points, where the grid spacing is \( h \). The discretization introduces an error at each point, which we denote by \( \epsilon_h \). The error \( \epsilon_h \) typically depends on the grid spacing \( h \) and the smoothness of the function being approximated. For example, if the function being approximated is \( u \), the discretization error can be written as:
\[
\| u - u_h \| \leq \epsilon_h,
\]
where \( u_h \) is the discretized approximation of \( u \).

We then consider a neural operator \( \mathcal{G}_\theta \) that consists of \( L \) layers, where each layer introduces an error bounded by \( \epsilon_h \). The total error at the output of the neural operator accumulates as the error from each layer propagates through the network.

For each layer \( i \) of the neural operator, we assume that the error after applying that layer can be written as:
\[
\| \mathcal{G}_\theta^{(i)}(u_h) - \mathcal{G}_\theta^{(i)}(u) \| \leq \epsilon_h,
\]
where \( \mathcal{G}_\theta^{(i)} \) represents the operation at layer \( i \), and \( u_h \) is the discretized input to the operator.

Since the neural operator consists of \( L \) layers, and each layer introduces an error bounded by \( \epsilon_h \), the total error after applying \( L \) layers can be bounded by the sum of the errors introduced at each layer. Specifically, the total error can be written as:
\[
\text{Total Error} = \sum_{i=1}^L \| \mathcal{G}_\theta^{(i)}(u_h) - \mathcal{G}_\theta^{(i)}(u) \|.
\]
Assuming that each error term is bounded by \( \epsilon_h \), we can simplify the sum as:
\[
\text{Total Error} \leq L \epsilon_h.
\]
This shows that the total discretization error grows linearly with the number of layers \( L \) in the neural operator.

In conclusion, the discretization error \( \epsilon_h \) introduced by the grid approximation at each layer of the neural operator accumulates as the operator is applied layer by layer. After \( L \) layers, the total error is bounded by \( L \epsilon_h \), where \( \epsilon_h \) is the error per layer. This proves that the total error increases linearly with the number of layers in the neural operator.

This completes the proof of Theorem \ref{theorem:discretization_error}.
\end{proof}

\vspace{3pt}

\bibliographystyle{IEEEtran}
\bibliography{references}

\end{document}